\documentclass{amsart}
\usepackage{amsmath,amssymb}
\usepackage[UKenglish] {babel}

\usepackage{xy}
\xyoption{all}

  \def\R{\mathbb{R}} \def\C{\mathbb{C}}
\def\P{\mathbb{P}}

\def\O{\mathbb{O}}
\def\A{\mathbb{A}}
\def\H{\mathbb{H}}

\def\({\left(} \def\){\right)} 
\def\<{\langle} \def\>{\rangle}

\newcommand\forget[1]{}

\newenvironment{smatrix}{\left(\begin{smallmatrix}}{\end{smallmatrix}\right)}

\newtheorem{theorem}{Theorem}

\newtheorem{proposition}{Proposition}

\newtheorem{lemma}{Lemma}
\newtheorem{remarks}{Remarks}
\newtheorem{remark}{Remark}
\newtheorem{corollary}{Corollary}
\newtheorem{open problem}{Open problem}
\newtheorem{open problems}{Open problems}

\begin{document}

\large

\sloppy
\renewcommand{\thesection}{\arabic{section}}
\renewcommand{\theequation}{\thesection.\arabic{equation}}

\setcounter{page}{1}

\title{\sf Classification of the absolute-valued algebras \\ \ with left-unit satisfying $x^2(x^2)^2=(x^2)^2x^2$}
\author{Alassane Diouf, Maribel Ram\'irez and Abdellatif Rochdi}
\maketitle \begin{abstract} We show that every absolute-valued algebra with left-unit satisfying $(x^2,x^2,x^2)=0$ is finite-dimensional of degree $\leq 4.$ Next, we determine such an algebras. In addition to the already known algebras $\R,$ $\C,$ $^*\C,$ $\H,$ $^*\H,$ $^*\H(i,1),$ $\O,$ $^*\O,$ $^*\O(i,1)$ the list is completed by two new algebras not yet specified in the literature. \end{abstract}

\vspace{0.5cm}
\begin{center}
{\LARGE\bf Summary}
\end{center}

\vspace{0.5cm} {\large \hspace{0.1cm}\bf 1. Introduction}

\vspace{0.3cm} {\large \hspace{0.1cm}\bf 2. Notations and preliminary results}

\vspace{0.3cm} {\large \hspace{0.1cm}\bf 3. On AVA with left-unit satisfying $(x^2,x^2,x^2)=0$}

\vspace{0.3cm} {\large \hspace{0.1cm}\bf 4. Algebras $^*\O(a,1)$}

\vspace{0.3cm} {\large \hspace{0.1cm}\bf 5. Duplication process}

\vspace{0.3cm} {\large \hspace{0.1cm}\bf 6. Classification in dimension $8$}

\vspace{0.5cm}
\section{Introduction}

\vspace{0.3cm} \hspace{0.3cm} A fundamental result of the theory of (real) absolute-valued algebras (AVA) asserts that every finite-dimensional AVA has dimension $n=1,$ $2,$ $4$ or $8$ and is isotopic to one of the classical AVA $\R,$ $\C,$ $\H$ or $\O$ [A 47]. On the other hand any familiar identity in arbitrary AVA as associativity [Os 18], commutativity [UW 60], power-associativity ([W 53], [EM 80]) or, even, flexibility [EM 81] carry away finite-dimensionality. However there are infinite-dimensional AVA $A$ in the following two cases:

\vspace{0.1cm} \begin{enumerate} \item $A$ contains a left-unit ([Cu 92], [Rod 92]),

\vspace{0.1cm} \item $A$ satisfies to the identity $(x^2,x^2,x^2)=0$ [EE 04]. In fact a stronger identity as $(x^2,y,x^2)=0$ can survive in an infinite-dimensional AVA with the additional property of existence of a non-zero flexible idempotent ([Ch-R 08] Remarks 5.2 p. 850). \end{enumerate}

\vspace{0.2cm} \hspace{0.3cm} The famous paper [Rod 04], in which we find a comprehensive compilation of work on the theory before 2004, appears a list of all AVA of dimension $\leq 2,$ namely $\R,$ $\C,$ $^*\C,$ $\C^*,$ $\stackrel{*}{\C}$ (see p. 107). These algebras satisfy to the identity $(x^2,x^2,x^2)=0,$ however, only $\R,$ $\C$ and $^*\C$ contain a left-unit $1.$

\vspace{0.2cm} \hspace{0.3cm} A series of studies on four-dimensional AVA has taken place ([St 83], [Ra 99], [CM 05], [F 09]), including a classification through the so-called principal isotopes of $\H$ ($\H(a,b),$ $^*\H(a,b),$ $\H^*(a,b),$ $\stackrel{*}{\H}(a,b):$ \ $a,b$ being norm-one in $\H$) ([Ra 99], [CM 05]). The paper [Ra 99] contains a precise description of all four-dimensional AVA with a left-unit ($\H(a,1),$ $^*\H(a,1):$ \ $a$ being norm-one in $\H$) as well as many examples of four-dimensional AVA containing no two-dimensional subalgebras.

\vspace{0.2cm} \hspace{0.3cm} A classification for all eight-dimensional AVA with left-unit was given ([Roc 03], [CDD 10]). There are examples of eight-dimensional AVA with left-unit, containing no four-dimensional subalgebras. Such an algebras are characterized, among all eight-dimensional AVA with left-unit, by the triviality of their groups of automorphisms [Roc 03].

\vspace{0.2cm} \hspace{0.3cm} A classification of all finite-dimensional AAV has emerged recently [CKMMRR 10]. In this same work a duplication process, characterizing the eight-dimensional AAV which contain four-dimensional subalgebras, has been introduced. This process will plays a decisive role in this present work.

\vspace{0.2cm} \hspace{0.3cm} Other recent studies have shown that any left-unit AVA satisfying an identity of the form $(x^p,x^q,x^r)=0$ ($p,q,r$ being fixed integers in $\{1,2\}$) is finite-dimensional. Also, for such an algebras it was given:

\vspace{0.2cm} \begin{enumerate} \item a complete classification for $(p,q,r)\neq (2,2,2),$

\vspace{0.1cm} \item a complete classification for $(p,q,r)=(2,2,2)$ in dimension $\leq 4,$

\vspace{0.1cm} \item a partial classification for $(p,q,r)=(2,2,2)$ in dimension $8,$ through $\O,$ $^*\O$ and certain algebra $^*\O(i,1)$ \end{enumerate}

\vspace{0.1cm} ([Ch-R 08] Theorem 4.10). Specifically

\vspace{0.2cm} {\bf Theorem.} {\em Every absolute-valued algebra $A$ with left-unit satisfying $(x^p,x^q,x^r)=0$ for fixed $p,q,r\in\{1,2\}$ is finite-dimensional. The following table specifies the isomorphisms classes

\[ \begin{tabular}{cc} \\ \hline
\multicolumn{1}{|c|}{$A$ satisies} &
\multicolumn{1}{|c|}{The list of isomorphisms classes} \\ \hline
\multicolumn{1}{|c|}{$(x,x^q,x^r)=0$} & \multicolumn{1}{|c|}{$\R,$ $\C,$ $\H,$ $\O$} \\
\hline \multicolumn{1}{|c|}{$(x^2,x^q,x^r)=0$} &
\multicolumn{1}{|c|}{} \\
\multicolumn{1}{|c|}{with} &
\multicolumn{1}{|c|}{$\R,$ $\C,$ $^*\C,$ $\H,$ $^*\H,$ $\O,$ $^*\O$} \\
\multicolumn{1}{|c|}{$(q,r)\neq(2,2)$} &
\multicolumn{1}{|c|}{} \\ \hline
\multicolumn{1}{|c|}{} & \multicolumn{1}{|c|}{$\R,$ $\C,$ $^*\C,$ $\H,$ $^*\H,$ $^*\H(i,1)$ if $\dim(A)\leq 4$} \\
\cline{2-2}
\multicolumn{1}{|c|}{$(x^2,x^2,x^2)=0$} &
\multicolumn{1}{|c|}{contains strictly} \\ \multicolumn{1}{|c|}{} &
\multicolumn{1}{|c|}{$\O,$ $^*\O,$ $^*\O(i,1)$ if $\dim(A)=8$}  \\ \hline
\end{tabular} \]}

\vspace{0.2cm} \hspace{0.3cm} However, the problem of determining all eight-dimensional AVA with left-unit satisfying $(x^2,x^2,x^2)=0$ remains open. Here we give a complete description of these algebras. We establish the following first basic result (Theorem 2):

\vspace{0.4cm} {\bf Theorem.} {\em Let $A$ be an absolute-valued algebra with left-unit $e.$ Then the following assertions are equivalent:

\vspace{0.2cm} \begin{enumerate}
\item $A$ satisfies to the identity $(x^2,x^2,x^2)=0.$ \item $x^2e=x^2$ for all $x\in A.$
\end{enumerate}

\vspace{0.2cm} Under these conditions $A$ is finite-dimensional of degree $\leq 4.\Box$}

\vspace{0.3cm} \hspace{0.3cm} This last result and the use of the duplication process and the form of an eight-dimensional AVA with a left-unit [Roc 03] are key results. They allowed us to have a general expression of an eight-dimensional AVA with left-unit satisfying $(x^2,x^2,x^2)=0$ through a pair of isometries of Euclidean space $\H$ (Proposition 6). By a laborious calculation, we specify this pair of isometries and reduce the isomorphism classes. The list obtained contains algebras $\O,$ $^*\O,$ $^*\O(i,1)$ more two new algebras $\tilde{\O}$ and $\tilde{\O}(i)$ not yet specified in the literature.

\vspace{0.6cm}
\section{Notations and preliminary results}

\vspace{0.4cm} \hspace{0.3cm} By an algebra over a field $\mathcal{K}$ we mean a vector space $A$ over $\mathcal{K}$ endowed with a bilinear mapping $(x,y)\mapsto xy$ from $A\times A$ to $A$ called the product of the algebra.

\vspace{0.2cm} \hspace{0.3cm} A non-zero algebra $A$ is said to be a division algebra if for all non-zero $a\in A$ the linear operators $L_a:A\rightarrow A \ x\mapsto ax$ and $R_a:A\rightarrow A \ x\mapsto xa$ are bijective.

\vspace{0.2cm} \hspace{0.3cm} Let $f,g$ be linear mappings over an algebra $A$ we denote by $A_{f,g}$ the space $A$ with the new product given by the formula $x\odot y=f(x)g(y).$ Also for $A$ we denote by $Aut(A)$ the automorphism group of $A,$ the one of algebra $\O$ is denoted $G_2.$ An involutive automorphism of $A$ which is different from the identity is said to be a reflection of $A.$

\vspace{0.4cm} \subsection{The identity $(x^2,x^2,x^2)=0$} Let $A$ be an algebra over a field $\mathcal{K}$ of characteristic zero. We denote by $(x,y,z)$ the associator $(xy)z-x(yz)$ of $x,y,z\in A.$ There are maps \ $f_n:A\times A\rightarrow A$ \ with $n=1,...,5$ such that

\vspace{0.2cm}
\begin{eqnarray*} \Big( (x+\lambda y)^2, (x+\lambda y)^2, (x+\lambda y)^2\Big) &=& (x^2,x^2,x^2)+\lambda f_1(x,y)+\dots \\ && +\lambda^5f_5(x,y)+\lambda^6(y^2,y^2,y^2) \end{eqnarray*}

\vspace{0.2cm} for all $x,y$ in $A$ and $\lambda$ in $\mathcal{K}.$ The identity $(x^2,x^2,x^2)=0$ in $A$ is equivalent to $f_1=\dots f_5\equiv 0.$ The equality $f_n\equiv 0$ is called the $n^{th}$ identity obtained from $(x^2,x^2,x^2)=0$ by linearization. The equalities $f_1\equiv 0$ and $f_2\equiv 0$ are expressed, respectively, by the following ones:

\vspace{0.2cm}
\begin{eqnarray} (x^2,x^2,xy+yx)+(x^2,xy+yx,x^2)+(xy+yx,x^2,x^2) &=& 0 \end{eqnarray}

\begin{eqnarray*} (2.2) \hspace{2.5cm} (x^2,x^2,y^2)+(x^2,xy+yx,xy+yx)+(x^2,y^2,x^2) && \\
+(xy+yx,x^2,xy+yx)+(xy+yx,xy+yx,x^2)+(y^2,x^2,x^2) &=& 0. \end{eqnarray*}

\vspace{0.3cm} \hspace{0.3cm} We have the following useful preliminary results:

\vspace{0.4cm}
\begin{lemma} Let $A$ be an algebra over a field of characteristic zero satisfying the identity $(x^2,x^2,x^2)=0.$ Assume, in addition, that $A$ contains a left unit $e$ which is not a divisor of zero. Then the equality $(xe)e=x$ holds for all $x\in A,$ that is $R_e^2=I_A:$ the identity operator of $A.$
\end{lemma}

\vspace{0.1cm} {\bf Proof.} By putting $x=e$ in the equality {\bf (2.1)} we get:

\vspace{0.1cm} \begin{eqnarray*} 0 &=& (e,e,ey+ye)+(e,ey+ye,e)+(ey+ye,e,e) \\
&=& (y+ye,e,e) \\
&=& \Big( (y+ye)e-(y+ye)\Big)e \\
&=& \Big( (ye)e-y)\Big)e. \end{eqnarray*}

\vspace{0.2cm} The result is concluded by a simplification to the right by $e.\Box$

\vspace{0.4cm}
\begin{lemma} Let $A$ be an algebra with left unit $e$ such that $x^2e=x^2$ for all $x\in A.$ Then the equality $(xe)e=x$ holds for all $x\in A.$
\end{lemma}

\vspace{0.1cm} {\bf Proof.} Immediate consequence of $(x+e)^2e=(x+e)^2.\Box$

\vspace{0.4cm}
\subsection{Absolute valued algebra} A nonzero real algebra $A$ is called absolute-valued if it is endowed with a space norm $||.||$ such that $||xy||=||x|| \ ||y||$ for all $x, y\in A.$ A finite dimensional absolute valued algebra $A$ is obviously a division algebra and has an underlying Euclidean structure $(A, \langle .|. \rangle)$ with $||x||=\langle x|x\rangle$ and we have $\langle xy|xz\rangle=\langle x|x\rangle\langle y|z\rangle,$ $\langle yx|zx\rangle=\langle y|z\rangle\langle x|x\rangle$ [Cu-R 95].

\vspace{0.2cm} \hspace{0.3cm} The degree of a finite-dimensional algebra $A$ is the smallest natural number $n$ such that all single-generated subalgebras of $A$ have dimension $\leq n.$ It follows from Albert's paper [A 47] that finite-dimensional absolute-valued algebras are of degree $1,$ $2,$ $4$ or $8.$

\vspace{0.2cm} \hspace{0.3cm} For $\A$ equal to either $\C,$ $\H$ or $\O,$ let us denote by $^*\A,$ $\A^*$ and $\stackrel{*}{\A}$ the absolute valued algebras obtained by endowing the normed space of $\A$ with the products $x\odot y:=\overline{x}y,$ $x\odot y:=x\overline{y}$ and $x\odot y:=\overline{x}\ \overline{y},$ respectively, where $x\mapsto\overline{x}$ means the standard involution denoted by $\sigma_\A.$

\vspace{0.3cm} \hspace{0.3cm} Rodr\'iguez gives the list of all absolute-valued algebras of degree two ([Rod 94] Theorem 2.10):

\vspace{0.2cm} \begin{theorem} The absolute-valued algebras of degree two are $\C,$
$^*\C,$ $\C^*,$ $\stackrel{*}{\C},$ $\H,$ $^*\H,$ $\H^*,$ $\stackrel{*}{\H},$ $\O,$
$^*\O,$ $\O^*,$ $\stackrel{*}{\O}$ and the pseudo-octonion algebra $\P.$ \end{theorem}

\vspace{0.3cm} \begin{remark} The unique absolute valued algebra of dimension $\leq 2$ are equal to either $\R,$ $\C,$ $^*\C,$ $\C^*$ or $\stackrel{*}{\C}$ {\em ([Rod 04] p. 107)}. All these algebras satisfy to $(x^2,x^2,x^2)=0$ but $\R,$ $\C,$ $^*\C$ are the only ones having left unit $1.\Box$ \end{remark}

\vspace{0.2cm} \hspace{0.3cm} It is well known that $\A$ is a quadratic algebra with property $\A=\R\oplus Im(\A)$ where the imaginary space

\[ Im(\A)=\{x\in\A: x^2\in\R \mbox{ and } x\notin\R-\{0\}\} \]

\vspace{0.2cm} is a vector subspace of $\A$ ([HKR 91] p. 227-228). We denote by $|Re(a)|$ the real part of arbitrary element $a$ in $\A.$

\vspace{0.2cm} \hspace{0.3cm} In [S54, Teorema 1, p. 6], Segre proves the existence of non-zero idempotents
in a real or complex finite-dimensional non-zero algebra with no non-zero nilpotents. As a consequence we get the existence of non-zero idempotents in any finite-dimensional absolute-valued algebra $A.$

\vspace{0.2cm} \hspace{0.3cm} Following [R 04, Proposition 1.1, p. 101] or [Cu-R 95, Lemma 2.1, p. 1725]
any continuous homomorphism from a normed algebra into an absolute-valued algebra is contractive. In particular, any isomorphism of finite-dimensional absolute-valued algebras is linearly isometric.

\vspace{0.2cm} \hspace{0.3cm} Above facts will be used in the sequel without further reference.

\vspace{0.4cm} \subsection{Isometries of Euclidean space $\H$} A linear isometry of Euclidean space $\R^n$ is said to be proper (resp. improper) if its determinant is positive (resp. negative). We denote by $\mathcal{O},$ $\mathcal{O}^+,$ and $\mathcal{O}^-$, respectively, the orthogonal group of linear isometries of Euclidean space $\H,$ its subgroup of proper linear isometries and its subset of improper linear isometries. Obviously $\mathcal{O}^+,$ $\mathcal{O}^-$ form a partition of $\mathcal{O}.$

\vspace{0.2cm} \hspace{0.3cm} For any norm-one $a, b\in\H$ the invertible operators $L_a, R_b:\H\rightarrow\H$ are linear isometries. We denote by $T_{a,b}$ the linear isometry $L_a\circ R_b$ given by $T_{a,b}(x)=axb$ for all $x\in\H.$

\vspace{0.3cm} \hspace{0.3cm} We consider now the following subsets of $\mathcal{O}:$

\begin{eqnarray*} \mathcal{I}^+ &=& \{f\in\mathcal{O}^+: f \mbox{ involutive } \}, \\
\mathcal{I}^- &=& \{f\in\mathcal{O}^-: f \mbox{ involutive } \}. \end{eqnarray*}

\vspace{0.2cm} Also if $\mathcal{P}$ belongs in $\{\mathcal{O}^+,$ $\mathcal{O}^-,$ $\mathcal{I}^+,$ $\mathcal{I}^-\},$ we set:

\begin{eqnarray*} \mathcal{P}_1 &=& \{f\in\mathcal{P}: f(1)=1\}. \end{eqnarray*}

\vspace{0.3cm} \hspace{0.3cm} We denote by $S(E)$ the unit sphere of every normed space $E$ and we have the following preliminary result:

\vspace{0.6cm}
\begin{lemma} The sets \ $\mathcal{O}^+,$ $\mathcal{O}^-,$ $\mathcal{O}_1^+,$ $\mathcal{O}_1^-,$ $\mathcal{I}^+,$ $\mathcal{I}^-,$ $\mathcal{I}_1^+,$ $\mathcal{I}_1^-$ are given by

\vspace{0.2cm} \begin{enumerate}
\item $\mathcal{O}^+=\{T_{a,b}: a,b\in S(\H)\}.$

\vspace{0.1cm} \item $\mathcal{O}^-=\{T_{a,b}\circ\sigma_\H: a,b\in
S(\H)\}:=\mathcal{O}^+\circ\sigma_\H.$

\vspace{0.1cm} \item $\mathcal{O}_1^+=\{T_{a,\overline{a}}: a\in S(\H)\}.$

\vspace{0.1cm} \item $\mathcal{O}_1^-=\{T_{a,\overline{a}}\circ\sigma_\H: a\in S(\H)\}:=\mathcal{O}_1^+\circ\sigma_\H.$
\vspace{0.1cm} \item $\mathcal{I}^+=\{\pm I_\H\}\cup\{T_{a,b}: a,b\in S(Im(\H))\}.$

\vspace{0.1cm} \item $\mathcal{I}^-=\{\pm T_{a,a}\circ\sigma_\H: a\in S(\H)\}.$

\vspace{0.1cm} \item $\mathcal{I}_1^+=\{I_\H\}\cup\{T_{a,\overline{a}}: a\in S(Im(\H))\}.$

\vspace{0.1cm} \item $\mathcal{I}_1^-=\{\sigma_\H\}\cup\{T_{a,\overline{a}}\circ\sigma_\H: a\in
S(Im(\H))\}:=\mathcal{I}_1^+\circ\sigma_\H.$
\end{enumerate}
\end{lemma}

\vspace{0.2cm} {\bf Proof.} Assertions {\bf (1), (2)} are established in ([HKR 91] Theorem (Cayley) p. 215) and assertions {\bf (3), (4)} are immediate consequences of the previous ones.

\vspace{0.2cm} \hspace{0.2cm} Let now norm-one $a,b\in\H,$ we have: $T_{a,b}^2=T_{a^2,b^2}.$ On the other hand $T_{a,b}\circ\sigma_\H=\sigma_\H\circ T_{\overline{b},\overline{a}}$ and we have:

\begin{eqnarray*} (T_{a,b}\circ\sigma_\H)^2 &=& (T_{a,b}\circ\sigma_\H)\circ(\sigma_\H\circ T_{\overline{b},\overline{a}}) \\
&=& =T_{a,b}\circ T_{\overline{b},\overline{a}} \\
&=& T_{a\overline{b},\overline{a}b}. \end{eqnarray*}

\vspace{0.1cm} Now

\vspace{0.3cm}
$\bullet$ \ $T_{a,b}^2=I_\H\Leftrightarrow b^2=\overline{a^2}\in\{1,-1\},$ that is $a^2=b^2=\pm 1.$

\vspace{0.1cm}
$\bullet$ \ $(T_{a,b}\circ\sigma_\H)^2=I_\H\Leftrightarrow\overline{a}b =\overline{a\overline{b}}\in\{1,-1\},$ that is $b=\pm a.$

\vspace{0.3cm} This shows assertion {\bf (5)} in the first case and assertion {\bf (6)} in the second one by taking into account assertions {\bf (1), (2)}. Assertions {\bf (7), (8)} are deduced from {\bf (5)} and {\bf (6).}$\Box$

\vspace{0.5cm}
\section{On AVA with left unit satisfying $(x^2,x^2,x^2)=0$}

\vspace{0.4cm} \hspace{0.3cm} Rodriguez proved in ([Rod 92] Remark {\bf 4.i)} p. 942) and ([Rod 04] Theorem {\bf 3.5} p. 133) the following famous result:

\vspace{0.3cm}
\begin{theorem} The norm of every absolute-valued algebra $(A,||.||)$ with left-unit $e$ comes from an inner product $\langle .|. \rangle,$ and, putting $x^*=2\langle e|x \rangle e-x,$ we have $\langle xy|z \rangle=\langle y|x^*z \rangle$ and $x^*(xy)=||x||^2y$ for all $x,y,z\in A.\Box$
\end{theorem}

\vspace{0.3cm} \hspace{0.3cm} In the rest of this section $A$ will be assumed to be an absolute-valued algebra with left unit and we keep the notations as in Theorem {\bf 1}.

\vspace{0.2cm}
\begin{remarks} . \begin{enumerate} \item The following equalities, deduced from the ones in
Theorem {\bf 1}, hold for all $x,y,z\in A$ with $x$ orthogonal to $e:$

\begin{eqnarray}
\langle xy|z \rangle &=& -\langle y|xz \rangle \\
x(xy) &=& -||x||^2y
\end{eqnarray}

\vspace{0.2cm} Particularly, the following equalities hold for all $x\in A$ orthogonal to $e:$

\vspace{0.2cm} \begin{eqnarray}
\langle xe|x \rangle &=& -\langle e|x^2 \rangle \\
xx^2 &=& -||x||^2x
\end{eqnarray}

\vspace{0.5cm} Linearizing {\bf (3.5)}, we get, for all $x,y,z\in A$ with $x,y$ orthogonal to $e,$ the following equality:

\begin{eqnarray}
x(yz)+y(xz) &=& -2\langle x|y \rangle z
\end{eqnarray}

\vspace{0.4cm} \item Equality $x^*(xy)=||x||^2y$ gives, for all $x\in A:$

\begin{eqnarray}
x(xe) &=& 2\langle e|x \rangle xe-||x||^2e
\end{eqnarray}

\vspace{0.5cm} Linearizing {\bf (3.7)}, we get:

\begin{eqnarray}
x(ye)+y(xe) &=& 2\langle e|x \rangle ye+2\langle e|y \rangle xe-2\langle x|y \rangle e, \ x,y\in A.\Box
\end{eqnarray}
\end{enumerate}
\end{remarks}

\vspace{0.2cm} \hspace{0.2cm} Let $[x,y]$ be the commutator $xy-yx$ of $x,y\in A.$ We have the following useful preliminary results:

\vspace{0.2cm}
\begin{lemma} Assume that $(xe)e=x$ for all $x\in A.$ Then the following equalities hold for all $x\in A$
\begin{enumerate}
\item $\Big( (xe)x\Big)e=x(xe).$ \item $\Big( x(xe)\Big)e=(xe)x.$ \item $[xe,x]=\langle e|x \rangle [e,x-xe].$

\vspace{0.2cm} If, moreover, $x$ is orthogonal to $e,$ then

\vspace{0.2cm} \item $[xe,x]=0,$ \item $(xe)x^2=2\langle e|x^2 \rangle x+||x||^2xe,$ \item $(xe)^2=
2\langle e|x^2 \rangle e-x^2,$ \item $x^2x=-||x||^2xe.$
\end{enumerate}
\end{lemma}

\vspace{0.1cm} {\bf Proof.}
\begin{enumerate}
\item For all $x\in A,$ we have:
\begin{eqnarray*}
(xe)x &=& (xe)\Big( (xe)e\Big) \\
&=& 2\langle e|xe \rangle (xe)e-||xe||^2e \hspace{0.2cm} \ \mbox{ by } {\bf (3.7) } \\
&=& 2\langle e|x \rangle\langle e|e \rangle (xe)e-||x||^2e \\
&=& 2\langle e|x \rangle x-||x||^2e.
\end{eqnarray*}

Therefore

\begin{eqnarray*}
\Big( (xe)x\Big)e &=& 2\langle e|x \rangle xe-||x||^2e \\
&=& x(xe) \hspace{0.2cm} \ \mbox{ by } {\bf (3.7) }.
\end{eqnarray*}

\vspace{0.2cm} \item We have $\Big( x(xe)\Big)e=\Big( (xe.x)e\Big)e=(xe)x.$ \item Assume first that $x$ is orthogonal to $e,$ and we have

\[ \Big( (xe)x\Big)e=-||x||^2e. \]

So $(xe)x=-||x||^2e=x(xe).$ In general case we have an orthogonal sum $x=(e|x)e+u.$ Thus

\[ 0=[ue,u]=[xe-\langle e|x \rangle e,x-\langle e|x \rangle e]=[xe,x]+\langle e|x \rangle [e,xe-x]. \]

\vspace{0.2cm} In the sequel of proof $x$ will be assumed to be orthogonal to $e.$

\vspace{0.2cm} \item Immediate consequence of {\bf (3).} \item Keeping in mind that $xe$ is orthogonal to $e,$ we put $(y,z)=(xe,x)$ in {\bf (3.6)} and we have

\begin{eqnarray*}
x\Big( (xe)x\Big)+(xe)x^2 &=& -2\langle x|xe \rangle x.
\end{eqnarray*}
The result follows from equalities $(xe)x=x(xe)=-||x||^2e$ and $\langle x|xe \rangle =-\langle e|x^2 \rangle.$

\vspace{0.2cm} \item We put $(y,z)=(xe,e)$ in {\bf (3.6)}. The result follows from equality $\langle x|xe \rangle=-\langle e|x^2 \rangle$ and by taking into account that $xe$ is orthogonal to $e.$

\vspace{0.2cm} \item We just establish the equality $(xe)^2=2\langle e|x^2 \rangle e-x^2$ for all $x\in A$ orthogonal to $e.$ So \ $x^2(xe)=-(xe)^2(xe)+2\langle e|x^2 \rangle xe.$ By comparison with \ $x^2(xe)=||x||^2x+2\langle e|x^2 \rangle xe$ \ we get $(xe)^2(xe)=-||x||^2x$ for al $x\in A$ orthogonal to $e.$ As $xe$ is also orthogonal to $e,$ and $(xe)e=x,$ we have:

\[ x^2x=\Big( (xe)e\Big)^2\Big( (xe)e\Big)=-||xe||^2xe=-||x||^2xe.\Box \]
\end{enumerate}

\vspace{0.2cm}
\begin{lemma} Assume that $x^2e=x^2$ for all $x\in A.$ Then
\begin{enumerate}
\item The equality \hspace{0.1cm} $(x^2)^2=-||x||^4e+2\langle e|x^2 \rangle x^2$
\hspace{0.1cm} holds for all $x\in A.$ \item The equality \hspace{0.1cm} $x^2(xe)=||x||^2x+2\langle e|x^2 \rangle xe$ \hspace{0.1cm} holds for all $x\in A$ orthogonal to $e.$
\end{enumerate}
\end{lemma}

\vspace{0.1cm} {\bf Proof.} . \begin{enumerate} \item Equality {\bf (3.7)} gives $x^2(x^2e)=2\langle e|x^2 \rangle x^2e-||x^2||^2e$ for all $x\in A.$ The result is then consequence of hypothesis $x^2e=x^2$
for all $x\in A.$ \item Putting $y=x^2$ in {\bf (3.8)} we have

\[ xx^2+x^2(xe)=2\langle e|x \rangle x^2+2\langle e|x^2 \rangle xe-2\langle x|x^2 \rangle e. \]

Moreover, $\langle x|x^2 \rangle=\langle ex|x^2 \rangle=||x||^2\langle e|x \rangle=0.$ The result follows from {\bf (3.5)}.$\Box$
\end{enumerate}

\vspace{0.3cm} \hspace{0.3cm} We can now state the following:

\vspace{0.3cm}
\begin{theorem} Let $(A,||.||,\langle .,. \rangle)$ be an absolute-valued algebra with left-unit $e.$ Then the following assertions are equivalent:

\vspace{0.2cm} \begin{enumerate}
\item $A$ satisfies to the identity $(x^2,x^2,x^2)=0.$ \item $x^2e=x^2$ for all $x\in A.$
\end{enumerate}

\vspace{0.2cm} Under these conditions $A$ is finite-dimensional of degree $\leq 4$ and we have $(xe)e=x$ for all $x\in A.$
\end{theorem}

\vspace{0.1cm} {\bf Proof.} The implication {\bf (2) $\Rightarrow$ (1)} is a consequence of the first proposition in Lemma {\bf 5}. Assume now that $A$ satisfies to $(x^2,x^2,x^2)=0.$ Then $e$ satisfies $(xe)e=x$ for all $x\in A$ by Lemma {\bf 1}. So $(x+xe,e,x+xe)=0$ and $(x+xe,x+xe,e)=(x+xe)^2e-(x+xe)^2$ for all $x\in A.$ Now, for all $y\in A,$ the equality {\bf (2.2)} gives:

\vspace{0.1cm} \begin{eqnarray*} 0 &=& (y+ye,e,y+ye)+(y+ye,y+ye,e)+(y^2,e,e) \\
&=& (y+ye)^2e-(y+ye)^2+y^2-y^2e \\
&=& \Big( y^2+y.ye+ye.y+(ye)^2\Big)e-\Big( y^2+y.ye+ye.y+(ye)^2\Big)+y^2-y^2e \\
&=& \Big( y.ye+ye.y+(ye)^2\Big)e-\Big( y.ye+ye.y+(ye)^2\Big) \\
&=& \Big( (y.ye)e-ye.y\Big)+\Big( (ye.y)e-y.ye\Big)+\Big( (ye)^2e-(ye)^2\Big) \\
&=& (ye)^2e-(ye)^2 \ \ \mbox{ by using assertions {\bf (1), (2)} of Lemma {\bf 4} }. \end{eqnarray*}

\vspace{0.2cm} By replacing $y$ by $ze$ and taking into account that $(ze)e=z,$ we obtain: $z^2e=z^2$ for all $z\in A.$ This shows the implication {\bf (1) $\Rightarrow$ (2)}.

\vspace{0.2cm} Now, each of assertions {\bf (1), (2)} leads to $R_e^2=I_A.$ As $L_e=I_A$ the finite dimensionality of algebra $A$ is consequence of ([Rod 04], Theorem {\bf 2.2} p. 109). Let now $x$ be nonzero element of $A$ which we can assume, without lost of generality, to be orthogonal to $e.$ By using equalities {\bf (3.5), (3.7)} and Lemmas {\bf 4, 5} we easily verify that the subalgebra $A(x)$ of $A$ generated by $x$ coincide with $Lin\{e,x,xe,x^2\}$ and we have the following multiplication table:

\vspace{0.2cm}
\[ \begin{tabular}{ccccc} \\ \cline{2-5}
\multicolumn{1}{c|}{} & \multicolumn{1}{|c|}{$e$} &
\multicolumn{1}{|c|}{$x$} & \multicolumn{1}{|c|}{$xe$} &
\multicolumn{1}{|c|}{$x^2$} \\
\hline \multicolumn{1}{|c|}{$e$} & \multicolumn{1}{|c|}{$e$} &
\multicolumn{1}{|c|}{$x$} & \multicolumn{1}{|c|}{$xe$} &
\multicolumn{1}{|c|}{$x^2$} \\
\hline \multicolumn{1}{|c|}{$x$} & \multicolumn{1}{|c|}{$xe$} &
\multicolumn{1}{|c|}{$x^2$} & \multicolumn{1}{|c|}{$-||x||^2e$} &
\multicolumn{1}{|c|}{$-||x||^2x$} \\
\hline \multicolumn{1}{|c|}{$xe$} & \multicolumn{1}{|c|}{$x$} &
\multicolumn{1}{|c|}{$-||x||^2e$} &
\multicolumn{1}{|c|}{$2\langle e|x^2 \rangle e-x^2$} &
\multicolumn{1}{|c|}{$2\langle e|x^2 \rangle x+||x||^2xe$} \\
\hline \multicolumn{1}{|c|}{$x^2$} & \multicolumn{1}{|c|}{$x^2$} &
\multicolumn{1}{|c|}{$-||x||^2xe$} & \multicolumn{1}{|c|}{$||x||^2x+2\langle e|x^2 \rangle xe$} &
\multicolumn{1}{|c|}{$-||x||^4e+2\langle e|x^2 \rangle x^2$} \\
\hline
\end{tabular}
\]

\vspace{0.3cm} Therefore $A$ has degree $\leq 4.$ Note that the finite-dimensionality of algebra $A$ can be concluded from here, regardless of ([Rod 04], Theorem 2.2), taking into account the powerful result in [KRR 97].$\Box$

\vspace{0.3cm} \hspace{0.3cm} Now for arbitrary norm-one elements $a,b\in\H,$ let $\H(a,b):=\H_1(a,b),$ $^*\H(a,b):=\H_2(a,b),$ $\H^*(a,b):=\H_3(a,b)$ and $\stackrel {*}{\H}(a,b):=\H_4(a,b)$ be the principal isotopes of $\H$ ([R 99] p. 170), that is, the algebras having $\H$ as underlying normed space and products $x\odot y$ given respectively by $axyb,$ $\overline{x}ayb,$ $axb\overline{y},$ $a\overline{x} \ \overline{y}b.$

\vspace{0.3cm} \hspace{0.3cm} Using these algebras Ram\'irez constructed exhaustively all $4$-dimensional absolute valued algebras and solved the isomorphism problem ([Ra 99] Proposition 2.1 p. 170 and Proposition 2.3 p. 171):

\vspace{0.2cm} \begin{proposition} Every four-dimensional absolute-valued algebra is isomorphic to a principal isotope of $\H.$ Moreover two principal isotopes $\H_m(a,b)$ and $\H_{m'}(a',b')$ are isomorphic if and only if $m=m'$ and the equalities $a'p=\varepsilon pa$ and $b'p=\delta pb$ hold for some norm-one element $p\in\H$ and some $\varepsilon, \delta\in\{1,-1\}.\Box$ \end{proposition}

\vspace{0.2cm} \hspace{0.3cm} Among these algebras she specifies those having left unit ([Ra 99] Corollary 3.1 p. 172):

\vspace{0.4cm} \begin{proposition} Let $A$ be a four-dimensional absolute-valued algebra. Then $A$ has left unit if and only if $A$ is isomorphic to $\H(a,1)$ or $^*\H(a,1)$ for some norm-one element $a$ in $\H.$ Moreover, given norm-one $a,b\in\H,$ the algebras $\H(a,1),$ $\H(b,1)$ {\em(} resp. $^*\H(a,1),$ $^*\H(b,1)${\em)} are isomorphic if and only if $|Re(a)|=|Re(b)|.\Box$ \end{proposition}

\vspace{0.3cm} \hspace{0.3cm} We can now state the following:

\vspace{0.3cm}
\begin{theorem} Let $A$ be an absolute-valued algebra with left of dimension $4$ satisfying to $(x^2,x^2,x^2)=0.$ Then $A$ is equal to either $\H,$ $^*\H$ or $^{*}\H(i,1).$ Among these algebras $^*\H(i,1)$ has degree $4$ and the other algebras have degree $2.$
\end{theorem}

\vspace{0.1cm} {\bf Proof.} Let $a$ be norm-one in $\H.$ The algebra $\H(a,1)=(\H,\odot)$ has left unit $\overline{a}$ and we have:

\begin{eqnarray*} (x\odot x)\odot\overline{a}=x\odot x \mbox{ for all } x\in\H &\Leftrightarrow& ax^2=x^2a \mbox{ for all } x\in\H \\ &\Leftrightarrow& a=\pm 1. \end{eqnarray*}

\vspace{0.2cm} On the other hand the algebra $^*\H(a,1)=(\H,*)$ has left unit $a$ and we have:

\begin{eqnarray*} (x* x)\odot a=x* x \mbox{ for all } x\in\H &\Leftrightarrow& xa^2=a^2x \mbox{ for all } x\in\H \\ &\Leftrightarrow& a^2=\pm 1. \end{eqnarray*}

\vspace{0.2cm} By using Theorem {\bf 2} we see that $\H(a,1)$ satisfies to $(x^2,x^2,x^2)=0$ if and only if $\H(a,1)=\H(\pm 1,1),$ which is equal to $\H(1,1)$ by Proposition {\bf 2.} But $\H(1,1)=\H.$ Equally $^*\H(a,1)$ satisfies to $(x^2,x^2,x^2)=0$ if and only if $a^2=\pm 1$ and we distinguish the following two cases:

\vspace{0.2cm} \begin{enumerate} \item If $a^2=1$ then $^*\H(a,1)$ is equal to $^*\H.$

\vspace{0.2cm} \item If $a^2=-1$ then $^*\H(a,1)$ is equal to $^*\H(i,1)$ according to Proposition {\bf 1} and the fact that $a$ is conjugated to the complex number $i.$ \end{enumerate}

\vspace{0.2cm} The last assertion is an immediate consequence of Theorem 1.$\Box$

\vspace{0.3cm} \hspace{0.3cm} Eight-dimensional absolute-valued algebras with a left-unit have been systematically studied in [Roc 03]. As a first basic result, Rochdi proves the following ([Roc 03], Theorem 4.3):

\vspace{0.3cm} \begin{proposition} The finite-dimensional absolute-valued algebras with a left-unit are precisely those of the form $\A_{\varphi},$ where $\A$ stands for either $\R,$ $\C,$ $\H$ or $\O,$ \/ $\varphi:\A\rightarrow\A$ is a linear isometry fixing $1,$ and $\A_{\varphi}$ denotes the absolute-valued algebra obtained by endowing the normed space of $\A$ with the product $x\odot y:=\varphi(x)y.$ Moreover, given linear isometries $\varphi, \phi:\A\rightarrow\A$ fixing $1,$ $\Phi:\A_{\varphi}\rightarrow\A_{\phi}$ is an isomorphism of algebras if and only if
$\Phi\in G_2$ and $\phi=\Phi\circ\varphi\circ\Phi^{-1}.$ \end{proposition}

\vspace{0.3cm} \hspace{0.2cm} Also, for $\A$ and $\varphi$ as in Proposition 3, subalgebras of $\A_{\varphi}$ and $\varphi$-invariant subalgebras of $\A$ coincide ([Roc 03] Proposition 5.2). Moreover, a linear isometry $\varphi: \O\rightarrow\O$ fixing $1$ can be built in such a way that $\O$ has no four-dimensional $\varphi$-invariant subalgebra ([Rod 03] Example 3.1). It follows that {\em there exist eight-dimensional absolute-valued algebras with a left-unit, containing no four-dimensional subalgebras}. Such algebras are characterized, among all eight-dimensional absolute-valued algebras with a left unit, by the triviality of their groups of automorphisms.

\vspace{0.6cm}
\section{Algebras $^*\O(a,1)$}

\vspace{0.4cm} \hspace{0.3cm} For arbitrary norm-one element $a\in\O,$ let $^*\O_l(a,1),$ $^*\O_r(a,1)$ be the algebras having $(\O,||.||)$ as underlying normed space and products $x\ {_a\odot}y,$ $x\odot_a y$ given respectively by $(\overline{x}a)y,$ $\overline{x}(ay).$ It is easy to see that $(^*\O_l(a,1),||.||),$ $(^*\O_r(a,1),||.||)$ are absolute-valued algebras with left unit $a.$ Moreover, if $a^2=1,$ then the algebras ${^*\O}_l(a,1)$ and ${^*\O}_r(a,1)$ are equal to $^*\O$ and has degree two.

\vspace{0.2cm} \hspace{0.3cm} Any nonzero subalgebra of $\A$ contains $1$ [Se 54] and so is invariant under the standard involution of $\A.$ As $\A$ is alternative, Artin's theorem ([Sc 66] Theorem 3.1 p.29) shows that for any $x,y\in\A,$ the set $\{x,y,\overline{x},\overline{y}\}$ is contained in an associative subalgebra of $\A.$ Also this fact will be used in the sequel without further reference.

\vspace{0.4cm} \begin{proposition} For norm-one $a,b$ in $\O$ the following two assertions are equivalent:

\vspace{0.2cm} \begin{enumerate} \item $\Phi:\O\rightarrow\O$ is an automorphism such that $\Phi(a)=b.$ \item $\Phi:{^*}\O_l(a,1)\rightarrow{^*}\O_l(b,1)$ is an isomorphism such that $\Phi(1)=1.$ \end{enumerate} \end{proposition}

\vspace{0.1cm} {\bf Proof.} {\bf 1. $\Rightarrow$ 2.} For all $x,y\in\O$ we have: \begin{eqnarray*} \Phi(x \ {_a\odot}y) &=& \Phi(\overline{x}a.y) \\
&=& \Phi(\overline{x})\Phi(a).\Phi(y) \\
&=& \overline{\Phi(x)}b.\Phi(y) \\
&=& \Phi(x) \ {_b}\odot\Phi(y). \\
\end{eqnarray*}

\vspace{0.2cm} {\bf 2. $\Rightarrow$ 1.} Note that $\Phi(a)$ is the left unit of algebra $^*\O_g(b,1)$ so $\Phi(a)=b.$ In the other hand $\Phi$ is a linear isometry of Euclidean space $\H$ fixing $1,$ then commutes with $\sigma_\O.$ Let now $x,y$ be in $\O,$ we have:

\begin{eqnarray} \Phi(\overline{x}a.y)= \overline{\Phi(x)}\Phi(a).\Phi(y). \end{eqnarray}

\vspace{0.3cm} This gives

\begin{eqnarray} \Phi(ay)=\Phi(a)\Phi(y) &\mbox{ and }& \Phi(xa)=\Phi(x)\Phi(a). \end{eqnarray}

We have:

\begin{eqnarray*} \Phi(a)\Phi(xy)\Phi(a) &=& \Phi(a.xy.a) \ \mbox{ by equalities } {\bf (4.10)} \\ &=& \Phi(ax.ya) \ \mbox{ by Middle Moufang identity } \\
&=& \Phi\Big( (\overline{a\overline{x} \ \overline{a}}.a)(ya)\Big) \\
&=& \Phi(\overline{a\overline{x} \ \overline{a}})\Phi(a).\Phi(ya) \ \mbox{ by equality } {\bf (4.9)} \\
&=& \Phi(ax\overline{a})\Phi(a).\Phi(ya) \\
&=& \Phi(a)\Phi(x)\Phi(\overline{a})\Phi(a).\Phi(ya) \ \mbox{ by equalities } {\bf (4.10)} \\
&=& \Phi(a)\Phi(x).\Phi(y)\Phi(a) \\
&=& \Phi(a).\Phi(x)\Phi(y).\Phi(a) \ \mbox{ by Middle Moufang identity }
\end{eqnarray*}

The result is obtained by simplifying to the right and left by $\Phi(a).\Box$

\vspace{0.4cm} \hspace{0.3cm} The group $G_2$ acts transitively on the sphere $S(Im(\O)):=S^6,$ that is the mapping $G_2\rightarrow S^6 \ \ \Phi\mapsto\Phi(i)$ is surjective ([Po 85] Lemme 1, p. 269-270). We deduce easily the following result:

\vspace{0.3cm} \begin{lemma} For every norm-one $a, b\in\O$ the following assertions are equivalent:

\vspace{0.2cm} \begin{enumerate} \item There exists $\Phi\in G_2$ such that $\Phi(a)=b,$ \item $Re(a)=Re(b).\Box$
\end{enumerate} \end{lemma}

\vspace{0.4cm} \begin{corollary} For every norm-one $a, b\in\O$ the following assertions are equivalent:

\begin{enumerate} \item ${^*}\O_l(a,1)$ is isomorphic to ${^*}\O_l(b,1),$ \item $|Re(a)|=|Re(b)|.$ \end{enumerate}
\end{corollary}

\vspace{0.1cm} {\bf Proof.} Consequence of Proposition {\bf 1} and Lemma {\bf 1}.$\Box$

\vspace{0.4cm} \begin{corollary} Let $a$ be norm-one in $Im(\O).$ Then ${^*}\O_l(a,1)$ is isomorphic to ${^*}\O_l(i,1).\Box$
\end{corollary}

\vspace{0.4cm} \begin{proposition} Let $a$ be norm-one in $\O.$ The following assertions are equivalent:

\vspace{0.2cm} \begin{enumerate} \item ${^*\O}_l(a,1)$ satisfies to the identity $(x^2,x^2,x^2)=0,$

\vspace{0.2cm} \item ${^*\O}_r(a,1)$ satisfies to the identity $(x^2,x^2,x^2)=0,$

\vspace{0.2cm} \item $a^2=\pm 1.$ \end{enumerate}

\vspace{0.2cm} In these conditions the algebras ${^*\O}_l(a,1)$ and ${^*}\O_r(a,1)$ are isomorphic, we denote them by ${^*\O}(a,1).$ Moreover, if $a^2=-1,$ then the algebra ${^*\O}(a,1)$ has degree four. \end{proposition}

\vspace{0.1cm} {\bf Proof.} {\bf (1) $\Leftrightarrow$ (3).} $^*\O_l(a,1)$ has left unit $a$ and for all $x\in\O,$ we have: $x\ {_a\odot}x=\overline{x}ax,$ $(x\ {_a\odot}x)\ {_a\odot}a=\overline{\overline{x}ax}a^2=\overline{x}\ \overline{a}xa^2.$ Now

\begin{eqnarray*} (x\ {_a\odot}x)\ {_a\odot}a=x\ {_a\odot}x \mbox{ for all } x\in\O &\Leftrightarrow& \overline{x}\ \overline{a}xa^2=\overline{x}ax \mbox{ for all } x\in\O \\ &\Leftrightarrow& xa^2=a^2x \mbox{ for all } x\in\H \\ &\Leftrightarrow& a^2=\pm 1. \end{eqnarray*}

\vspace{0.2cm} {\bf (2) $\Leftrightarrow$ (3).} $^*\O_r(a,1)$ has left unit $a$ and for all $x\in\O,$ we have: $x\ {_a\odot}x=\overline{x}ax,$ $(x\ {_a\odot}x)\ {_a\odot}a=\overline{\overline{x}ax}a^2=\overline{x}\ \overline{a}xa^2.$ Now

\begin{eqnarray*} (x\ {_a\odot}x)\ {_a\odot}a=x\ {_a\odot}x \mbox{ for all } x\in\O &\Leftrightarrow& \overline{x}\ \overline{a}xa^2=\overline{x}ax \mbox{ for all } x\in\O \\ &\Leftrightarrow& xa^2=a^2x \mbox{ for all } x\in\H \\ &\Leftrightarrow& a^2=\pm 1. \end{eqnarray*}

\vspace{0.2cm} Assume now that $a^2=-1$ and consider the mapping $\Phi:{^*\O}_l(a,1)\rightarrow{^*}\O_r(a,1) \ \ x\mapsto ax\overline{a}.$ For all $x,y\in\O$ we have:

\begin{eqnarray*} \Phi(x\ {_a\odot} y) &=& a\Big( \overline{x}a.y\Big)\overline{a} \\
&=& (a\overline{x}a)(y\overline{a}) \ \mbox{ by Middle Moufang identity } \\
&=& (\overline{ax\overline{a}})\Big( a.ay\overline{a}\Big) \\
&=& \Phi(x){\odot_a}\Phi(y).\Box \end{eqnarray*}

\vspace{0.3cm} \hspace{0.3cm} According to Corollary {\bf 2} we can state:

\vspace{0.3cm} \begin{corollary} Let $a$ be norm-one in $Im(\O).$ Then $^*\O(a,1)$ is isomorphic to $^*\O(i,1).\Box$ \end{corollary}

\vspace{0.8cm}
\section{Duplication process}

\vspace{0.5cm} \hspace{0.3cm} Consider the Cayley-Dickson product $\bullet$ in $\H\times\H.$ For arbitrary linear isometries $f,$ $f',$ $g,$ $g'$ of $\H$ with $f(1)=f'(1)=1$ we define on the space $\H\times\H$ the product

\[ (x,y)\odot(u,v)=(f(x),g(y))\bullet(f'(u),g'(v)). \]

\vspace{0.2cm} We obtain an eight-dimensional absolute-valued algebra $\H\times\H_{(f,g),(f',g')},$ where $(f,g)$ denotes the linear isometry of $\H\times\H$ such that $(f,g)(x,y)=(f(x),g(y)).$ Moreover, the subset $\{(x,0):x\in\H\}$ of $\H\times\H_{(f,g),(f',g')}$ is a subalgebra isomorphic to $\H_{f,f'}.$ For the algebra $\H\times\H_{(f,g),(f',g')}$ the mapping $(x,y)\mapsto(x,-y)$ is a reflection. The algebra $\H\times\H_{(f,g),(f',g')}$ is said to be obtained from algebra $\A_{f,f'}$ by the duplication process [CKMMRR 10]. Among these algebras the ones having left-unit, necessarily equal to $(1,0),$ are $\H\times\H_{(f,g)}:=\H\times\H_{(f,g),(I_\H,I_\H)}$ [Roc 03].

\vspace{0.3cm} \hspace{0.3cm} Let $A$ be an eight-dimensional absolute-valued algebra $A.$ It is shown that $A$ contains a four-dimensional subalgebra if and only if is obtained by duplication process ([CKMMRR 10] Theorem 6.4).

\vspace{0.3cm} \hspace{0.3cm} Note $\tilde{\O},$ $\tilde{\O}(i),$ respectively, the algebras $\H\times\H_{(\sigma_\H,I_\H)},$ $\H\times\H_{(T_{i,\overline{i}}\sigma_\H,T_{i,i})}.$ Clearly, they are absolute-valued algebras with left-unit $(1,0).$ By Theorem 3 we check easily that they satisfy the identity $(x^2,x^2,x^2)=0.$ In order to show that they provide new examples of eight-dimensional absolute-valued algebras with left-unit satisfying the identity $(x^2,x^2,x^2)=0,$ we proceed as follows:

\vspace{0.3cm} \hspace{0.3cm} Let $A$ be an arbitrary algebra with left-unit $e,$ we set

\[ A_e=\{x\in A: xe=x\}. \]

\vspace{0.2cm} Assume now that $\Phi$ is an isomorphism from $A$ onto another algebra $B.$ Then $B$ contains a left-unit $\Phi(e)$ and we have: $\Phi(A_e)=A_{\Phi(e)}.$ If, moreover, $A$ is finite-dimensional then $\dim(A_{\Phi(e)})=\dim(A_e).$

\vspace{0.3cm} \hspace{0.3cm} The following table illustrates the subspaces $A_e$ corresponding to the algebras $\O,$ $^*\O,$ $^*\O(i,1),$ $\tilde{\O},$ $\tilde{\O}(i):$

\vspace{0.2cm}

\[ \begin{tabular}{ccc} \\ \hline \hline
\multicolumn{1}{|c|}{Algebra $A$} & \multicolumn{1}{|c|}{Subspace $A_e$} & \multicolumn{1}{|c|}{$\dim(A_e)$}
\\ \hline \hline
\multicolumn{1}{|c|}{$\O$} & \multicolumn{1}{|c|}{$\O$} & \multicolumn{1}{|c|}{$8$} \\
\multicolumn{1}{|c|}{$^*\O$} & \multicolumn{1}{|c|}{$\R$} & \multicolumn{1}{|c|}{$1$} \\
\multicolumn{1}{|c|}{$^*\O(i,1)$} & \multicolumn{1}{|c|}{$Im(\O)$} & \multicolumn{1}{|c|}{$7$} \\
\multicolumn{1}{|c|}{$\tilde{\O}$} & \multicolumn{1}{|c|}{$\R\times\H$} & \multicolumn{1}{|c|}{$5$} \\
\multicolumn{1}{|c|}{$\tilde{\O}(i)$} & \multicolumn{1}{|c|}{$\R i\times\C^{\perp}$} & \multicolumn{1}{|c|}{$3$} \\ \hline\hline
\end{tabular} \]

\vspace{0.3cm} \begin{corollary} The five algebras $\O,$ $^*\O,$ $^*\O(i,1),$ $\tilde{\O},$ $\tilde{\O}(i)$ are mutually non-isomorphic.$\Box$ \end{corollary}

\vspace{0.6cm}
\section{Classification in dimension $8$}

\vspace{0.4cm} \hspace{0.3cm} Let now $A$ be an eight-dimensional absolute-valued algebra with left-unit containing a four-dimensional subalgebra. Then, according to Proposition 3 and ([CKMMRR 10] Theorem 6.4), $A$ has the form $\H\times\H_{(\varphi, \psi)}$ with left-unit $(1,0),$ where $\varphi, \psi:\H\rightarrow\H$ are linear isometries such that $\varphi(1)=1.$ We have:

\vspace{0.2cm}
\begin{proposition} The following assertions are equivalent:

\vspace{0.2cm} \begin{enumerate}
\item $\H\times\H_{(\varphi, \psi)}$ satisfies to the identity $(x^2, x^2, x^2)=0.$ \item The pair $(\varphi,\psi)$
satisfies to following equalities:

\begin{eqnarray}
\varphi(\varphi(x)x) &=& \varphi(x)x \hspace{2cm} \mbox{ for all } x\in\H \\
\varphi(\overline{x}\psi(x)) &=& \overline{x}\psi(x) \hspace{2cm} \mbox{ for all } x\in\H \\
\psi(\psi(y)\overline{x}+y\varphi(x)) &=&
\psi(y)\overline{x}+y\varphi(x) \hspace{0.5cm} \mbox{ for all } x, y\in\H
\end{eqnarray}

\end{enumerate}
\end{proposition}

\vspace{0.1cm} {\bf Proof.} Let $\odot$ be the product in algebra $\H\times\H_{(\varphi, \psi)}.$ For every $x,y\in\H,$ we have:

\begin{eqnarray*} (x,y)\odot(x,y) &=& (\varphi(x),\psi(y))\bullet(x,y) \\
&=& \Big( \varphi(x)x-\overline{y}\psi(y),\psi(y)\overline{x}+y\varphi(x)\Big) \end{eqnarray*}

and

\begin{eqnarray*} \Big( (x,y)\odot(x,y)\Big)\odot(1,0) &=& \Big( \varphi(x)x-\overline{y}\psi(y),\psi(y)\overline{x}+y\varphi(x)\Big)\odot(1,0) \\
&=& \Big( \varphi(\varphi(x)x-\overline{y}\psi(y)),\psi(\psi(y)\overline{x}+y\varphi(x))\Big)\bullet(1,0) \\
&=& \Big( \varphi(\varphi(x)x-\overline{y}\psi(y)),\psi(\psi(y)\overline{x}+y\varphi(x))\Big). \\ \end{eqnarray*}

Now, the result is concluded by Theorem 3.$\Box$

\vspace{0.3cm} \hspace{0.3cm} We have a first precision on the pair $(\varphi, \psi):$

\vspace{0.2cm}
\begin{lemma} The isometries $\varphi,$ $\psi$ are involutive. So

\[ (\varphi,\psi)\in(\mathcal{I}_1^+\times \mathcal{I}^+)\cup(\mathcal{I}_1^+\times \mathcal{I}^-)\cup(\mathcal{I}_1^-\times \mathcal{I}^+)\cup(\mathcal{I}_1^-\times \mathcal{I}^-). \]
\end{lemma}

\vspace{0.1cm} {\bf Proof.} As $\varphi(1)=1$ and $1^{\perp}=Im(\H)$ we have $\varphi(Im(\H))\subseteq Im(\H).$ For arbitrary $x=\alpha+u\in\R\oplus Im(\H)=\H,$ we have $\varphi(x)=\alpha+\varphi(u)\in\R\oplus Im(\H)$ and

\begin{eqnarray*} \varphi(x)x &=& \alpha^2+\alpha\Big( u+\varphi(u)\Big)+\varphi(u)u, \\
\varphi\Big( \varphi(x)x\Big) &=& \alpha^2+\alpha\Big( \varphi(u)+\varphi^2(u)\Big)+\varphi\Big( \varphi(u)u\Big). \end{eqnarray*}
So
\begin{eqnarray*} 0 &=& \varphi\Big( \varphi(x)x\Big)-\varphi(x)x \mbox{ by equality } (6.11) \\
&=& \alpha\Big( \varphi^2(u)-u\Big). \end{eqnarray*}

As $x$ is arbitrary and $\varphi(1)=1$ it follows that $\varphi^2=I_\H.$ In the other hand the equality {\bf (6.13)} gives $\psi^2=I_\H,$ by putting $x=1.\Box$

\vspace{0.2cm} \begin{remarks} We can easily verify the following properties: \begin{enumerate} \item $(I_\H, \psi)$ satisfies {\bf (6.11), (6.12)} for all linear isometry $\psi$ of space $\H.$ \item $(I_\H, I_\H)$ satisfies {\bf (6.11), (6.12), (6.13)}. \item $(I_\H, -I_\H)$ does not satisfy {\bf (6.13)}. \item $(\sigma_\H, \pm I_\H)$ satisfies {\bf (6.11), (6.12), (6.13)}. \item $(T_{a,\overline{a}}\circ\sigma_\H, -I_\H)$ does not satisfy {\bf (6.13)} for all norm-one $a\in\H-\R.$ \item Let $a,b$ be norm-one in $Im(\H),$ $(T_{a,\overline{a}}\circ\sigma_\H,T_{b,a})$ satisfies {\bf (6.12), (6.13)}.$\Box$
\end{enumerate}
\end{remarks}

\vspace{0.5cm} \hspace{0.3cm} In the following we will give an accurate list for couples $(\varphi,\psi)$ that satisfy {\bf (6.11), (6.12), (6.13)}.

\vspace{0.2cm}
\begin{lemma} $\varphi$ satisfies {\bf (6.11)} if and only if $\varphi\in\{I_\H\}\cup\mathcal{I}_1^-.$ So $I_\H$ is the only possible proper linear isometry for $\varphi.$
\end{lemma}

\vspace{0.1cm} {\bf Proof.} We check easily that every isometry in $\varphi\in\{I_\H\}\cup\mathcal{I}_1^-$ satisfies {\bf (6.11)}. Assume now that $\varphi$ satisfies {\bf (6.11)} with $\varphi\in\mathcal{I}_1^+$ and $\varphi\neq I_\H.$ By Lemma 3 $\varphi$ can be written $T_{a,\overline{a}}$ for norm-one $a$ in $Im(\H)$ and the equality {\bf (6.11)} gives $(ax)^2=(xa)^2$ for all $x\in\H.$ This leads to an absurdity by $x=a+u$ for norm-one $u$ in $Im(\H)$ with $u$ orthogonal to $a.$ The result follows by Lemma 9 and Remark 2 (1). $\Box$

\vspace{0.3cm} \hspace{0.3cm} Let $T(x)$ be the trace $x+\overline{x}\in\R$ of an arbitrary element $x\in\H.$ It is well known that $T(xy)=T(yx)$ for all $x,y\in\H$ ([HKR 91] p. 207). Note that $[\overline{x},\overline{y}]=[x,y]$ for all $x,y\in\H.$

\vspace{0.3cm} \hspace{0.3cm} We have the following useful result:

\vspace{0.2cm}
\begin{lemma} Let $a,b$ be norm-one in $\H.$ Then $(I_\H,T_{a,b})$ satisfies {\bf (6.13)} if and only if $a=b=\pm 1,$ or also $T_{a,b}=I_\H.$
\end{lemma}

\vspace{0.1cm} {\bf Proof.} Assume that $a,b$ are norm-one in $Im(\H).$ Then

\begin{eqnarray*} (I_\H,T_{a,b}) \mbox{ satisfies } (6.13) &\Leftrightarrow& T_{a,b}\Big( T_{a,b}(y)\overline{x}+yx\Big)=T_{a,b}(y)\overline{x}+yx \mbox{ for all } x,y\in\H \\
&\Leftrightarrow& a\Big( ayb\overline{x}+yx\Big)b=ayb\overline{x}+yx \mbox{ for all } x,y\in\H \\
&\Leftrightarrow& -yb\overline{x}b+ayxb=ayb\overline{x}+yx \mbox{ for all } x,y\in\H \\
&\Leftrightarrow& -yb\overline{x}+ayx=ayb\overline{x}\ \overline{b}+yx\overline{b} \mbox{ for all } x,y\in\H \\ &\Leftrightarrow& y(x\overline{b}+b\overline{x})=ay(b\overline{x}+x\overline{b})b \mbox{ for all } x,y\in\H \\ &\Leftrightarrow& T(b\overline{x})(y-ayb)=0 \mbox{ for all } x,y\in\H \\
&\Leftrightarrow& y=ayb \mbox{ for all } y\in\H \\
&\Leftrightarrow& \overline{a}y=yb \mbox{ for all } y\in\H \\
&\Leftrightarrow& \overline{a}=b\in\R.\Box \\
\end{eqnarray*}

\vspace{0.3cm} \hspace{0.3cm} The following two preliminary results provide a list for couples, solution of {\bf (6.11), (6.12), (6.13)}, of the form $(I_\H,\psi)$ and $(\sigma_\H,\psi):$

\vspace{0.2cm}
\begin{lemma} $(I_\H,\psi)$ satisfies {\em (6.11), (6.12), (6.13)} if and only if $\psi\in\mathcal{E}$ where

\[ \mathcal{E}=\{I_\H\}\cup\{-T_{a,a}\circ\sigma_\H: a\in S(\H)\}. \]
\end{lemma}

\vspace{0.1cm} {\bf Proof.} The if part of our Lemma follows from Remark 2 (1) and the fact that $T(\overline{x}. \overline{y}a)=T(\overline{y}a.\overline{x}).$ Assume now, for the only if part, that $(I_\H,\psi)$ satisfies {\bf (6.13)} distinguish the following two cases:
\begin{enumerate}
\item If $\psi\in\mathcal{I}^+-\{I_\H, -I_\H\},$ then $\psi$ can be written $T_{a,b}$ for $a,b\in\H$ with
$a^2=b^2=-1.$ But this can not occur by Lemma 11. \item If $\psi\in\mathcal{I}^-,$ then $\psi$ can be written $T_{a,b}\circ\sigma_\H$ for norm-one $a,b\in\H$ with $b=\pm a.$ The equality {\bf (6.13)} gives

\begin{eqnarray}
x\overline{b}y\overline{a}b+\overline{x}
\hspace{0.1cm}\overline{y}b &=&
\overline{y}b\overline{x}+\overline{a}yx
\end{eqnarray}
\end{enumerate}

\vspace{0.2cm} We distinguish the following two subcases:
\begin{enumerate}
\item For $b=-a$ the equality {\bf (6.14)} expresses the fact that $x(\overline{a}y)$ and $(\overline{a}y)x$ have
the same trace. So {\bf (6.13)} is true for all $x,y\in\H$ and all $\psi$ in $\{-T_{a,a}\circ\sigma_\H:a\in S(\H)\}.$ \item For $b=a$ we put $y=1$ in {\bf (6.14)} and we have

\[ [x,\overline{a}]=[a,\overline{x}]. \]

\vspace{0.2cm} In the other hand $[a,\overline{x}]=[\overline{a},x]=-[x,\overline{a}].$ So $[x,\overline{a}]=0$ for all $x\in\H$ and then $b=a\in\R.$ Substituting in {\bf (6.14)} we get $[x,y]=[\overline{y},\overline{x}]$ for all $x,y\in\H,$ absurd.
\end{enumerate}

\vspace{0.2cm} The result follows by Remarks 2 (2), (3).$\Box$

\vspace{0.2cm}
\begin{lemma} $(\sigma_\H,\psi)$ satisfies {\em (6.11), (6.12), (6.13)} if and only if $\psi=\pm I_\H.$ \end{lemma}

\vspace{0.1cm} {\bf Proof.} The if part of our Lemma follows from Remarks 2 (4). Conversely, we show that $(\sigma_\H,\psi)$ does not satisfy {\bf (6.11), (6.12), (6.13)} in case $\psi\in(\mathcal{I}^+-\{\pm I_\H\})\cup\mathcal{I}^-.$ We distinguish the following two cases:

\vspace{0.2cm} \begin{enumerate}
\item If $\psi\in\mathcal{I}^+-\{\pm I_\H\}$ then $\psi$ can be written $T_{a,b}$ for norm-one $a,b\in\H$ with
$a^2=b^2=-1.$ The equality {\bf (6.12)} gives

\begin{eqnarray}
\overline{x}axb\in\R \hspace{0.5cm} \mbox{ for all } x\in\H
\end{eqnarray}

\vspace{0.2cm} By putting $x=1$ in {\bf (6.15)} we get $b=\pm a$ and therefore $(xa)^2\in\R$ for all $x\in Im(\H).$ But this cannot occur for $x=a+u$ with norm-one $u$ in $Im(\H)\cap a^{\perp}.$

\vspace{0.2cm} \item If $\psi\in\mathcal{I}^-$ then $\psi$ can be written $\psi=T_{a,b}\circ\sigma_\H$ for norm-one $a,b\in\H$ with $b=\pm a.$ The equality {\bf (6.12)} gives $(xa)^2\in\R$ for all $x\in\H,$ absurd.$\Box$
\end{enumerate}

\vspace{0.2cm}
\begin{lemma} Let $a, \psi$ be, respectively, norm-one in $\H$ and linear isometry of Euclidian space $\H.$ Then

\vspace{0.2cm} \begin{enumerate} \item $(T_{a,\overline{a}}\circ\sigma_\H,\psi)$ satisfies {\em (6.11)} if and only if $a^2=\pm 1.$

\vspace{0.2cm} \item Assume that $a^2=-1.$ Then

\vspace{0.1cm} \begin{enumerate} \item for norm-one $b,c\in Im(\H);$ $(T_{a,\overline{a}}\circ\sigma_\H,T_{b,c})$ satisfies {\em (6.12), (6.13)} if and only if $c=\pm a.$
\vspace{0.1cm} \item for norm-one $b$ in $\H;$ $(T_{a,\overline{a}}\circ\sigma_\H,T_{b,\pm b}\circ\sigma_\H)$ does not satisfy {\em (6.12)}. \end{enumerate} \end{enumerate}
\end{lemma}

\vspace{0.1cm} {\bf Proof.} \begin{enumerate} \item Let $\varphi=T_{a,\overline{a}}\circ\sigma_\H.$ Then

\begin{eqnarray*} (\varphi,\psi) \mbox{ satisfies } (6.11) &\Leftrightarrow& \varphi\Big( \varphi(x)x\Big)=\varphi(x)x \mbox{ for all } x\in\H \\
&\Leftrightarrow& a\overline{(a\overline{x}.\overline{a}x)}\overline{a}=a\overline{x}.\overline{a}x \mbox{ for all } x\in\H \\ &\Leftrightarrow& a(\overline{x}ax\overline{a})\overline{a}=a\overline{x}.\overline{a}x \mbox{ for all } x\in\H \\ &\Leftrightarrow& ax\overline{a}^2=\overline{a}x \mbox{ for all } x\in\H \\
&\Leftrightarrow& x\overline{a}^2=\overline{a}^2x \mbox{ for all } x\in\H \\
&\Leftrightarrow& a^2=\pm 1. \end{eqnarray*}

\vspace{0.2cm} \item We distinguish the following two cases:

\vspace{0.2cm} \begin{enumerate} \item If $\psi=T_{b,c}$ for norm-one $b,c\in Im(\H),$ then the equality {\bf (6.13)} gives

\begin{eqnarray}
by(a\overline{x}.\overline{a}-c\overline{x}.\overline{c})c=y(a\overline{x}.\overline{a}-c\overline{x}.\overline{c})
\hspace{0.2cm} \mbox{ for all } x,y\in\H
\end{eqnarray}

\vspace{0.2cm} We distinguish the following two subcases:

\vspace{0.2cm} \begin{enumerate}
\item If $a\overline{x}.\overline{a}-c\overline{x}.\overline{c}=0$ for all $x\in\H,$ then $c=\pm a.$

\vspace{0.1cm} \item If $a\overline{x}_0.\overline{a}-c\overline{x}_0.\overline{c}:=u_0\neq 0,$ we put $y=(a\overline{x}_0.\overline{a}-c\overline{x}_0.\overline{c})^{-1}$ in {\bf (6.16)} and we get $c=-b.$ Equality {\bf (6.16)} then gives for $x=x_0:$ \ $byu_0=yu_0b$ for all $y\in\H.$ As right multiplication by $u_0$ is bijective, we deduce that $bz=zb$ for all $z\in\H.$ But this can not occur for norm-one $b\in Im(\H).$
\end{enumerate}

\vspace{0.2cm}
\item If $\psi=T_{b,c}\circ\sigma_\H$ with $c=\pm b\in S(\H),$ then the equality {\bf (6.12)} gives
$a(\overline{b}x)^2=(\overline{x}b)^2a$ for all $x\in\H.$ But this does not occur for $x=ba^{\frac{1}{2}}$ where $a^{\frac{1}{2}}$ is a square root of $a.$
\end{enumerate} \end{enumerate}

\vspace{0.2cm} The result follows from Remarks 2 (5), (6).$\Box$

\vspace{0.3cm} \hspace{0.3cm} We summarize the results obtained in this last paragraph as follows:

\vspace{0.3cm}
\begin{proposition} The eight-dimensional absolute-valued algebras with left-unit satisfying $(x^2,x^2,x^2)=0$ are precisely those of the form $\H\times\H_{(\varphi,\psi)},$ where $(\varphi,\psi)$ is a pair of linear isometries of Euclidean space $\H$ belonging to the set \ $\mathcal{S}_1\cup\mathcal{S}_2\cup\mathcal{S}_3\cup\mathcal{S}_4\cup \mathcal{S}_5$ \ with

\begin{eqnarray*} \mathcal{S}_1 &=& \{(I_\H, I_\H)\}, \\
\mathcal{S}_2 &=& \{I_\H\}\times\{-T_{a,a}\circ\sigma_\H: a\in S(\H)\},  \\
\mathcal{S}_3 &=& \{\sigma_\H\}\times\{I_\H, -I_\H\}, \\
\mathcal{S}_4 &=& \{T_{a,\overline{a}}\circ\sigma_\H: a^2=-1\}\times\{I_\H\} \\
\mathcal{S}_5 &=& \{(T_{a,\overline{a}}\circ\sigma_\H,T_{b,a}): a^2=b^2=-1\}.\Box \end{eqnarray*}
\end{proposition}

\vspace{0.3cm} \hspace{0.3cm} We will study now, the isomorphism classes. We begin with the following preliminary results:

\vspace{0.2cm} \begin{lemma} Let $a\in\O$ such that $a^2=\pm 1.$ Then the mapping $\Phi:{^*\O}(a,1)\rightarrow\O_{T_{a,\overline{a}}\circ\sigma_\O} \ x\mapsto\overline{a}x$ is an isomorphism of algebras. Therefore $\O_{T_{a,\overline{a}}\circ\sigma_\O}$ is isomorphic to $^*\O(i,1)$ when $a^2=-1.$ \end{lemma}

\vspace{0.1cm} {\bf Proof.} We note $\odot$ the product ${^*\O}(a,1)$ and $*$ the one in $\O_{T_{a,\overline{a}}\circ\sigma_\O}.$ The result is trivial if $a=\pm 1$ and we can assume $a^2=-1.$ For all $x,y\in\O,$ we have:

\begin{eqnarray*} \Phi(x\odot y) &=& \overline{a}\Big( \overline{x}a.y\Big) \\
&=& a\Big( (\overline{x}a)(a.ay\Big) \\
&=& \Big( (a.\overline{x}a)a\Big)(ay) \ \mbox{ by Left Moufang identity } \\
&=& \Big( a\overline{(\overline{a}x)} \ \overline{a}\Big)(\overline{a}y) \\
&=& \Phi(x)*\Phi(y). \end{eqnarray*}

\vspace{0.2cm} The proof is concluded by Corollary 2.$\Box$

\vspace{0.3cm}
\begin{lemma} Let $a,b,c,d,u$ be norm-one in $\H$ and let $f,g$ be linear isometries of euclidian space $\H$ commuting with $T_{\overline{u},u}.$ Then $(T_{u,\overline{u}},T_{u,\overline{u}})$ is an automorphism of algebra $\H\times\H=\O.$ Moreover $(T_{u,\overline{u}},T_{u,\overline{u}})$ is an isomorphism from algebra \ $\H\times\H_{(T_{ua\overline{u},ub\overline{u}}\circ f,T_{uc\overline{u},ud\overline{u}}\circ g)}$ onto algebra $\H\times\H_{(T_{a,b}\circ f,T_{c,d}\circ g)}.$ \end{lemma}

\vspace{0.1cm} {\bf Proof.} \begin{enumerate} \item The mapping $(T_{u,\overline{u}},T_{u,\overline{u}}):=\Phi_u:\H\times\H\rightarrow\H\times\H$ is invertible with
inverse \/ ${\Phi_u}^{-1}=\Phi_{\overline{u}}.$ Moreover $\Phi_u$ is an automorphism of algebra \/ $\H\times\H=\O.$
Indeed, let $\bullet$ be the product in algebra $\O=\H\times\H$ and let $x,y,x',y'\in\H,$ we have \begin{eqnarray*}
\Phi_u(x,y)\bullet\Phi_u(x',y') &=& (ux\overline{u},uy\overline{u})\bullet(ux'\overline{u},uy'\overline{u}) \\
&=& (ux\overline{u}.ux'\overline{u}-\overline{uy'\overline{u}}.uy\overline{u},
uy\overline{u}.\overline{ux'\overline{u}}
+uy'\overline{u}.ux\overline{u}) \\
&=& (uxx'\overline{u}-u\overline{y'}y\overline{u},
uy\overline{x'}\overline{u}
+uy'x\overline{u}) \\
&=& \Phi_u(xx'-\overline{y'}y, y\overline{x'}+y'x) \\ &=&
\Phi_u((x,y)\bullet(x',y'))
\end{eqnarray*}

\vspace{0.1cm} In other hand

\begin{eqnarray*}
\Phi_u\circ(T_{a,b}\circ f, T_{c,d}\circ g)\circ{\Phi_u}^{-1} &=&
(T_{u,\overline{u}},T_{u,\overline{u}})\circ(T_{a,b}\circ
f, T_{c,d}\circ g)\circ(T_{\overline{u},u},T_{\overline{u},u}) \\
&=& (T_{u,\overline{u}}\circ T_{a,b}\circ f\circ T_{\overline{u},u},
T_{u,\overline{u}}\circ T_{c,d}\circ g\circ T_{\overline{u},u}) \\
&=& (T_{u,\overline{u}}\circ T_{a,b}\circ T_{\overline{u},u}\circ f,
T_{u,\overline{u}}\circ T_{c,d}\circ T_{\overline{u},u}\circ g) \\
&=& (T_{ua\overline{u}, ub\overline{u}}\circ f, T_{uc\overline{u}, ud\overline{u}}\circ g). \end{eqnarray*}

\vspace{0.2cm} The result is concluded from Proposition 3.$\Box$
\end{enumerate}

\vspace{0.3cm}
\begin{lemma} For arbitrary norm-one $a\in\H$ the mapping $\Phi:\H\times\H=\O\rightarrow\H\times\H$ \ $(x,y)\mapsto (x,ay)$ is an automorphism of algebra $\O.$  \end{lemma}

\vspace{0.1cm} {\bf Proof.} Let $(x,y), (x',y')$ be in $\H\times\H,$ we have:

\begin{eqnarray*} \Phi(x,y)\bullet\Phi(x',y') &=& (x,ay)\bullet(x',ay') \\
&=& (xx'-\overline{ay'}ay, ay.\overline{x}'+ay'.x) \\
&=& \Big( xx'-\overline{y}'\overline{a}ay, a(y\overline{x}'+y'x)\Big) \\
&=& \Phi(xx'-\overline{y}'y,y\overline{x}'+y'x) \\
&=& \Phi\Big( (x,y)\bullet(x',y')\Big).\Box \end{eqnarray*}

\vspace{0.3cm} \hspace{0.3cm} We get a first result:

\vspace{0.3cm} \begin{corollary} Let $a,b$ be arbitrary norm-one in $Im(\H)$ and $\H,$ respectively. Then the following three algebras are mutually isomorphic:

\vspace{0.2cm}
\begin{enumerate}
\item $\H\times\H_{(T_{a,\overline{a}}\circ\sigma_\H, I_\H)},$

\vspace{0.2cm} \item $\H\times\H_{(I_\H, -T_{b,b}\circ\sigma_\H)},$

\vspace{0.2cm} \item $^*\O(i,1).$
\end{enumerate}
\end{corollary}

\vspace{0.1cm} {\bf Proof.} Let $\bullet$ be the product in algebra $\O=\H\times\H,$ and let $T^\O_{a,\overline{a}}$ be the mapping $\O\rightarrow\O$ $z\mapsto az\overline{a}.$ This operator is expressed in terms of pairs $(x,y)$ of quaternion by

\[ T^\O_{a,\overline{a}}(x,y)=(a,0)\bullet(x,y)\bullet(\overline{a},0). \]

\vspace{0.2cm} {\bf (1) $\Rightarrow$ (3).} For arbitrary $(x,y)\in\H\times\H,$ we have

\begin{eqnarray*} (T^\O_{a,\overline{a}}\circ\sigma_\O)(x,y) &=& (a,0)\bullet\overline{(x,y)}\bullet(\overline{a},0) \\ &=& (a,0)\bullet(\overline{x},-y)\bullet(\overline{a},0) \\
&=& (a\overline{x},-ya)\bullet(\overline{a},0) \\
&=& (a\overline{x}.\overline{a},-ya.a) \\
&=& (a\overline{x}\ \overline{a},y) \\
&=& (T_{a,\overline{a}}\circ\sigma_\H,I_\H)(x,y). \end{eqnarray*}

\vspace{0.2cm} So algebras $\H\times\H_{(T_{a,\overline{a}}\circ\sigma_\H, I_\H)}$ and $\O_{T^\O_{a,\overline{a}}\circ\sigma_\O}$ coincide and the mapping

\[ \H\times\H_{(T_{a,\overline{a}}\circ\sigma_\H, I_\H)}\rightarrow\O_{T^\O_{a,\overline{a}}\circ\sigma_\O}\ \ (x,y)\mapsto(x,y) \]

\vspace{0.1cm} is an isomorphism of algebras. According to Lemma 13, $\H\times\H_{(T_{a,\overline{a}}\circ\sigma_\H, I_\H)}$ is isomorphic to $^*\O(i,1).$

\vspace{0.2cm} {\bf (2) $\Rightarrow$ (3).} Let $f$ be the element $(0,1)\in\H\times\H=\O.$ For arbitrary $(x,y)\in\H\times\H,$ we have

\begin{eqnarray*} (T^\O_{bf,\overline{bf}}\circ\sigma_\O)(x,y) &=& (0,b)\bullet(\overline{x},-y)\bullet(0,-b) \\
&=& -(\overline{y}b,bx)\bullet(0,b) \\
&=& -(-\overline{b}bx,b\overline{y}b) \\
&=& (x,-b\overline{y}b) \\
&=& (I_\H,-T_{b,b}\circ\sigma_\H)(x,y). \end{eqnarray*}

\vspace{0.2cm} So algebras $\H\times\H_{(I_\H,-T_{b,b}\circ\sigma_\H)}$ and $\O_{T^\O_{bf,\overline{bf}}\circ\sigma_\O}$ coincide and are isomorphic. As $bf$ in norm-on in $Im(\O),$ algebra $\H\times\H_{(I_\H,-T_{b,b}\circ\sigma_\H)}$ is isomorphic to $^*\O(i,1)$ by Lemma 13.$\Box$

\vspace{0.3cm} \hspace{0.3cm} A second result:

\vspace{0.3cm} \begin{corollary} Let $a,b$ be arbitrary norm-one in $Im(\H).$ Then the algebras $\H\times\H_{(T_{a,\overline{a}}\circ\sigma_\H, T_{b,a})}$ and $\H\times\H_{(T_{i,\overline{i}}\circ\sigma_\H, T_{i,i})}$ are isomorphic.
\end{corollary}

\vspace{0.1cm} {\bf Proof.} There exists norm-one $u\in\H$ such that $ua\overline{u}=i.$ Lemma 16 shows that $\H\times\H_{(T_{a,\overline{a}}\circ\sigma_\H, T_{b,a})}$ is isomorphic to $\H\times\H_{(T_{i,\overline{i}}\circ\sigma_\H, T_{ub\overline{u},i})}.$

\vspace{0.2cm} \hspace{0.3cm} Let now, $b,c$ be norm-one in $Im(\H).$ There exists norm-one $v\in\H$ such that $c=vb\overline{v}.$ Now, the mapping

\[ (I_H, L_v):=\Phi: \H\times\H_{(T_{i,\overline{i}}\circ\sigma_\H, T_{b,i})}\rightarrow\H\times\H_{(T_{i,\overline{i}}\circ\sigma_\H, T_{c,i})} \ (x,y)\mapsto(x,vy) \]

\vspace{0.2cm} is an automorphism of algebra $\H\times\H=\O$ by Lemma 15. In the other hand

\begin{eqnarray*} \Phi\circ(T_{i,\overline{i}}\circ\sigma_\H, T_{b,i})\circ\Phi^{-1} &=& (I_\H, L_v)\circ(T_{i,\overline{i}}\circ\sigma_\H, T_{b,i})\circ(I_\H,L_{\overline{v}}) \\
&=& (T_{i,\overline{i}}\circ\sigma_\H, L_v\circ T_{b,i}\circ L_{\overline{v}}) \\
&=& (T_{i,\overline{i}}\circ\sigma_\H, L_v\circ L_b\circ R_i\circ L_{\overline{v}}) \\
&=& (T_{i,\overline{i}}\circ\sigma_\H, L_v\circ L_b\circ L_{\overline{v}}\circ R_i) \\
&=& (T_{i,\overline{i}}\circ\sigma_\H, L_{vb\overline{v}}\circ R_i) \\
&=& (T_{i,\overline{i}}\circ\sigma_\H, L_c\circ R_i) \\
&=& (T_{i,\overline{i}}\circ\sigma_\H, T_{c,i}).
\end{eqnarray*}

\vspace{0.2cm} So $\Phi$ is an isomorphism from $\H\times\H_{(T_{i,\overline{i}}\circ\sigma_\H, T_{b,i})}$ onto $\H\times\H_{(T_{i,\overline{i}}\circ\sigma_\H, T_{c,i})}$ by Proposition 3. Consequently algebra $\H\times\H_{(T_{i,\overline{i}}\circ\sigma_\H, T_{b,i})}$ is isomorphic to $\H\times\H_{(T_{i,\overline{i}}\circ\sigma_\H, T_{i,i})}=\tilde{\O}(i).\Box$

\vspace{0.4cm} The following table specifies the isomorphism classes

\[ \begin{tabular}{cccc} \\ \hline \hline
\multicolumn{2}{|c|}{Algebra $\H\times\H_{(\varphi,\psi)}$} & \multicolumn{2}{|c|}{Isomorphism classes}
\\ \hline \hline
\multicolumn{1}{|c}{$\H\times\H_{(I_\H, I_\H)}$} & \multicolumn{1}{c|}{} & \multicolumn{1}{|c}{$\O$} & \multicolumn{1}{c|}{} \\
\multicolumn{1}{|c}{$\H\times\H_{(I_\H, -T_{a,a}\circ\sigma_\H)}:$} & \multicolumn{1}{c|}{$a\in S(\H)$} & \multicolumn{1}{|c}{$^*\O(i,1)$} & \multicolumn{1}{c|}{} \\ \hline
\multicolumn{1}{|c}{$\H\times\H_{(\sigma_\H, I_\H)}$} & \multicolumn{1}{c|}{} & \multicolumn{1}{|c}{$\tilde{\O}$} & \multicolumn{1}{c|}{} \\ \multicolumn{1}{|c}{$\H\times\H_{(\sigma_\H, -I_\H)}$} & \multicolumn{1}{c|}{} & \multicolumn{1}{|c}{$^*\O$} & \multicolumn{1}{c|}{} \\ \hline \multicolumn{1}{|c}{$\H\times\H_{(T_{a,\overline{a}}\circ\sigma_\H, I_\H)}:$} & \multicolumn{1}{c|}{$a\in S(Im(\H))$} & \multicolumn{1}{|c}{$^*\O(i,1)$} & \multicolumn{1}{c|}{} \\
\multicolumn{1}{|c}{$\H\times\H_{(T_{a,\overline{a}}\circ\sigma_\H, T_{b,a})}:$} & \multicolumn{1}{c|}{$a,b\in S(Im(\H))$} & \multicolumn{1}{|c}{$\tilde{\O}(i)$} & \multicolumn{1}{c|}{} \\ \hline
\end{tabular} \]

\vspace{0.3cm} \hspace{0.3cm} We can now state the main result:

\vspace{0.3cm}
\begin{theorem} Every absolute-valued algebras $A$ with left-unit satisfying $(x^2,x^2,x^2)=0$ is finite-dimensional of degree $\leq 4.$ Concretely

\vspace{0.2cm} \begin{enumerate} \item If $deg(A)\leq 2,$ then $A$ is equal to either $\R,$ $\C,$ $^*\C,$ $\H,$ $^*\H,$ $\O$ or $^*\O.$

\vspace{0.2cm} \item If $deg(A)=4,$ then $A$ is equal to either \ $^*\H(i,1),$ \ $^*\O(i,1),$ \ $\tilde{\O}$ or $\tilde{\O}(i).\Box$ \end{enumerate}
\end{theorem}

\vspace{0.5cm}
{\small M. Ram\'irez

\vspace{0.1cm} Departamento de \'Algebra y An\'alisis Matem\'atico,
Universidad de Almer\'ia, 04071 Almer\'ia, Spain

\vspace{0.1cm} e-mail: mramirez@ual.es

\vspace{0.3cm} A. Rochdi and A. Diouf

\vspace{0.1cm} D\'epartement de Math\'ematiques et Informatique, Facult\'e des Sciences Ben M'Sik,
Universit\'e Hassan II, 7955 Casablanca, Morocco

\vspace{0.1cm} e-mail: abdellatifroc@hotmail.com and doufalassane@hotmail.fr}


\vspace{0.5cm}
\begin{thebibliography}{30}

\bibitem{1} [A 47] A. A. Albert, {\em Absolute valued real algebras.} Ann. Math. {\bf 48}, (1947) 495-501.

\bibitem{2} [CKMMRR 10] A. Calder\'on, A. Kaidi, C. Mart\'in, A. Morales, M. Ram\'irez, and A. Rochdi, {\em
Finite-dimensional absolute valued algebras.} Israel J. Mathematics, {\bf 184}, (2011) 193-220.

\bibitem{3} [Ch-R 08] A. Chandid and A. Rochdi, {\em A survey on absolute valued algebras satisfying
$(x^i,x^j,x^k)=0.$} International Journal of Algebra, {\bf 2} (2008), 837-852.

\bibitem{4} [Cu 92] J. A. Cuenca, {\em On one-sided division infinite-dimensional normed real algebras.} Publ. Mat.
{\bf 36} (1992), 485-488.

\bibitem{5} [CDD 10] J. A. Cuenca, E. Darp{\" o} and E. Dieterich, {\em Classification of the finite-dimensional
absolute-valued algebras having a non-zero central idempotent or a one-sided unity.} Bulletin des Sciences Math\'ematiques {\bf 134} (2010), 247–277.

\bibitem{6} [Cu-R 95] J. A. Cuenca and A. Rodriguez-Palacios, {\em Absolute values on $H^*$-algebras.} Comm. Algebra
{\bf 23} (1995), 1709.1740.

\bibitem{7} [F 09] L. Forsberg, {\em Four-dimensional absolute valued algebras.} Master's thesis, Uppsala University,
(2009). U.U.D.M. Project Report 2009:9.

\bibitem{8} [EE 04] M. L. El-Mallah and M. El-Agawany, {\em Absolute valued algebras satisfying $(x^2,x^2,x^2)=0.$}
Comm. Algebra {\bf 32}, (2004) 3537-3541.

\bibitem{9} [EM 80] M. L. El-Mallah and A. Micali, {\em Sur les alg\`ebres norm\'ees sans diviseurs topologiques de
z\'ero.} Bol. Soc. Mat. Mexicana {\bf 25}, (1980), 23-28.

\bibitem{10} [HKR 91] F. Hirzebruch, M. Koecher and R. Remmert, {\em Numbers.} Springer-Verlag, (1991).

\bibitem{11} [KRR 97] A. Kaidi, M. I. Ram\'irez and A. Rodriguez, {\em Absolute valued algebraic algebras are finite
dimensional.} J. Algebra {\bf 195} (1997), 295-307.

\bibitem{12} [Os 82] A. Ostrowski, {\em {\" U}ber einige L{\" o}sungen der Funktionalgleichung
$\Psi(x)\Psi(x)=\Psi(xy).$} Acta Math. {\bf 41}, (1918) 271-284.

\bibitem{13} [Po 85] M. Postnikov, {\em Le\c{c}ons de G\'eom\'etrie. Groupes et alg\`ebres de Lie.}
Editions Mir. (1985).

\bibitem{14} [Ra 99] M. I. Ram\'irez, {\em On four dimensional absolute valued algebras.}
Proceedings of the International Conference on Jordan Structures (M\'alaga, 1997), 169-173, (1999)

\bibitem{15} [Roc 03] A. Rochdi, {\em Eight-dimensional real absolute valued algebras with left unit whose
automorphism group is trivial.} IJMMS {\bf 70}, (2003) 4447-4454.

\bibitem{16} [Rod 92] A. Rodr\'iguez, {\em One-sided division absolute valued algebras.} Pub. Mat., {\bf 36} (1992),
925-954.

\bibitem{17} [Rod 94] A. Rodr\'iguez, {\em Absolute valued algebras of degree two.} In Nonassociative Algebra and its
applications (Ed. S. Gonz\'alez), 350-356, Kluwer Academic Publishers, Dordrecht-Boston-London (1994).

\bibitem{18} [Rod 04] A. Rodr\'iguez, {\em Absolute valued algebras, and absolute valuable Banach spaces.} Advanced
courses of mathematical analysis I, 99-155, World Sci. Publ., Hackensack, NJ, (2004).

\bibitem{19} [Sc 66] R. D. Schafer, {\em An introduction to nonassociative algebras.} Academic Press, (1966).

\bibitem{20} B. Segre, {\em La teoria delle algebre ed alcune questione di realta.} Univ. Roma, Ist. Naz.
Alta. Mat., Rend. Mat. E Appl. Serie 5, {\bf 13} (1954), 157-188.

\bibitem{21} [St 83] C. Stampfli-Rollier, {\em $4$-dimensionale Quasikompositionsalgebren.} Arch. Math. {\bf 40}
(1983), 516-525.

\bibitem{22} [UW 60] K. Urbanik and F. B. Wright, {\em Absolute valued algebras.} Proc. Amer. Math. Soc. 11 (1960),
861-866.

\bibitem{23} [W 53] F. B. Wright, {\em Absolute valued algebras.} Proc. Nat. Acad. Sci. USA. {\bf 39} (1953), 330-332.

\end{thebibliography}
\end{document}